\documentstyle[12pt,amssymbols]{article}  
\newcommand{\bq}{\begin{equation}}  
\newcommand{\eq}{\end{equation}}  
\newcommand{\bqa}{\begin{eqnarray}}  
\newcommand{\eqa}{\end{eqnarray}}  
\newcommand{\ra}{\rightarrow}

\def\ed{\end{document}}  
  
\def\ra{\rightarrow}  
\def\al{\alpha}  
\def\2pi{1\over 2\pi i}  
\def\q{q-q^{-1}}  
\def\~{\tilde}
\def\newline{\hfil\break}  
  
\def\la{\lambda}  
\def\ra{\rightarrow}

\def\sq2{\sqrt{2}}  
\def\sqk2{\sqrt{2(k+2}}  
\def\sqk{\sqrt{k}}

\def\be{\begin{equation}}  
\def\ee{\end{equation}}  
\def\br{\begin{array}}  
\def\er{\end{array}}  
\def\bea{\begin{eqnarray}}  
\def\eea{\end{eqnarray}}  
\def\ba{\begin{equation}\begin{array}}  
\def\ea{\end{array}\end{equation}}  
\def\bac{\begin{equation}\begin{array}{rll}}

\def\al{\alpha}

\newcommand{\uq}{U_q (\widehat{sl(2)})}

\def\Z{{\Bbb Z}}  
\def\C{{\Bbb C}}

\begin{document}  
\rightline{CRM-2311}  
\rightline{ August, 1995}  
\vbox{\vspace{-10mm}}  
\vspace{1.0truecm}  
\begin{center}  
{\LARGE \bf  
On the relation between $U_q(\widehat{sl(2)})$  
vertex operators and $q$-zonal functions
 }\\[8mm]  
{\large A.H. Bougourzi\footnote{Address as of September 1, 
1995: Institute for Theoretical Physics, State University of 
New York at Stony Brook, Stony Brook, N.Y. 11794, USA} and L. Vinet  }\\  
[6mm]{\it 
Centre de Recherches Math\'ematiques\\  
Universit\'e de Montr\'eal\\  
C.P. 6128-A, Montr\'eal, P.Q., H3C 3J7, Canada.  
}\\[20mm]  
  
\end{center}  
\vspace{1.0truecm}  
\begin{abstract}      

We show how the states constructed from the action of the modes
 of bosonized vertex operators, that intertwine 
$U_q(\widehat{sl(2)})$ modules, are related to $q$-zonal
functions.
\end{abstract} 

\newpage

\section{Introduction}  
 
The space of states of an exactly-solvable two-dimensional model 
typically consists of a finite number of infinite-dimensional 
highest weight modules of an affine algebra or of its 
$q$-deformation. 
These modules 
can be constructed using either the currents, or the
intertwining vertex operators. They can also be embedded in 
the Fock
spaces constructed from  the action of the 
Heisenberg subalgebras and certain 
group algebras on the weight lattice of the  Lie 
subalgebra.

It has been shown \cite{BakerThesis} that the currents 
generating $\uq$ 
are related to certain symmetric functions introduced by Kerov 
in \cite{K}.
Since the correlation functions of exactly-solvable models 
are expectation values of the intertwiners, it is natural to 
relate the latter to symmetric functions also. 

Let $\Lambda_i$, $i=\{0,1\}$ be the $sl(2)$ fundamental weights, 
$V(\Lambda_i)$  the level-one $\uq$ highest weight modules, 
$V$  the $U_q(sl(2))$ two-dimensional module with basis vectors 
$v_\pm$, $z$  a complex variable, and 
$V(z)=V \otimes {\C}[z,z^{-1}]$ 
 the level-zero $\uq$-evaluation module. 
In \S 2,  we show that the states constructed from 
the modes of the vertex operators $\phi^{i,\pm}(z)$, 
 referred to as type-I in \cite{DFJMN} and which intertwine 
the  $\uq$-modules as   
\be  
\phi^{i,+}(z)\otimes v_{+}+\phi^{i,-}(z)\otimes v_{-}:\quad  
V(\Lambda_i)\ra V(\Lambda_{i-1})\otimes V(z), 
\label{boubi}\ee  
are related to the Macdonald functions 
$P_{\lambda}(x;q^{4},q^{2})$ \cite{MacdonaldBook}.
In the limit $q \ra 1$, these functions reduce to the zonal    
functions $Z_\lambda(x)$ and we shall hence refer to them 
as the  $q$-zonal   functions.  (The  functions 
$Z_\lambda(x)$ are related to the Jack polynomials 
[5] in the following way $Z_\lambda(x) = J_\lambda(x;2)$.) 
In \S 3, we compare our results with those of Jing \cite{Jing1} 
who considered some general vertex operators, related to 
Macdonald 
functions,  but {\em not} to $\uq$. In particular, 
Jing's operators are not related to intertwiners of $\uq$ 
modules. 
We conclude by discussing a few open problems. 
 
\section{Relation between $q$-zonal   functions and  
$\uq$ intertwining vertex operators} 
 
In this section, we review the bosonization of the  vertex 
operators (VO) (1.1), following \cite{JMMN}, and consider the 
states of the Fock space constructed by acting  with  the 
modes of these 
VO on the bosonic vacuum state. We refer to these  
states as {\it the vertex 
operator states (VOS)}.  The VO $\phi^{i,\pm}(z)$ are 
bosonized as follows 
\cite{JMMN}:  
\bac  
\phi^{i,-}(z)&=&  
e^{\sum_{n=1}^{\infty}{q^{7n/2} \over [2n]}a_{-n}z^{n}}  
e^{-\sum_{n=1}^{\infty}{q^{-5n/2}\over  
[2n]}a_{n}z^{-n}}e^{\al/2}(-q^3z)^{(\partial+i)/2}
=\sum_{n\in\Z} 
\phi^{i,-}_{n} z^{-n},\\  
\phi^{i,+}(z)&=&-{\q\over 2\pi i}\oint_{|zq^4|<|w|<|zq^2|} 
dw {w  
:\phi^{i,-}(z)  
E^-(w):\over   
(w-zq^2)(w-zq^4)}=\sum_{n\in\Z} 
\phi^{i,+}_{n} z^{-n},  
\ea  
where $[x]=(q^{x}-q^{-x})/(q-q^{-1})$. The normal ordering 
symbol 
$::$ means that the creation modes ${a_n, e^\al; n<0}$ are 
placed 
to the left of the annihilation modes ${a_n, \partial; n>0}$, 
and 
$E^-(z)$ is one of the quantum currents that generate $\uq$.   
The explicit realization of this current is given by  
\be  
E^-(z)=e^{-\sum_{n=1}^{\infty}{q^{n/2}\over [n]}a_{-n}z^{n}}  
e^{\sum_{n=1}^{\infty}{q^{n/2}\over  
[n]}a_{n}z^{-n}}e^{-\al}z^{-\partial},  
\ee  
where $\al$ is the simple positive root of $sl(2)$. Moreover, 
$e^\al$ is a translation operator acting on the $sl(2)$ weight 
lattice, while $\partial$ is the zero mode (momentum) operator 
which is 
conjugate to $\al$ in the following sense:  
\bac  
{[\partial,\al]}&=&2,\\  
z^{\partial}e^{\alpha}&=&z^2:e^{\alpha}z^{\partial}:.  
\ea  
In actual computations invoking to the zero mode $\partial$ and 
$\al$, only the last relation is used. The non-zero modes 
satisfy the following commutation relations:  
\be  
{[a_n,a_m]}={[2n][n]\over n}\delta_{n+m,0}.  
\ee  
For the purposes of \S 3, where we compare our results with
Jing's, it is useful to redefine the non-zero modes as follows:  
\bac  
a_n&\equiv& {[n]\over n}q^{-|n|}b_n,\quad 
\forall n\in\Z\backslash\{0\},\\  
{[b_n,b_m]}&=&n(1+q^{2|n|})\delta_{n+m,0}.  
\label{HA} 
\ea  
In terms of the modes $b_n$, 
the  vertex operators are then expressed as  follows 
\bac  
\phi^{i,-}(z)&=&  
e^{\sum_{n=1}^{\infty}{q^{4n}\over n(1+q^{2n})}b_{-n}z^{n}}  
e^{-\sum_{n=1}^{\infty}{q^{-2n}\over  
n(1+q^{2n})}b_{n}z^{-n}}e^{\al/2}(-q^3z)^{(\partial+i)/2}\\
&=& \sum_{n\in\Z}\phi^{i,-}_{n}  z^{-n},\\  
\phi^{i,+}(z)&=&-{\q\over 2\pi i}\oint_{|zq^4|<|w|<|zq^2|} dw 
{w  :\phi^{i,-}(z)  
E^-(w):\over   
(w-zq^2)(w-zq^4)
},  
\label{vo}\ea  
with   
\bac  
:\phi^{i,-}(z)E^-(w):=&e&^{\sum_{n=1}^{\infty}{q^{4n}\over 
n(1+q^{2n})}b_{-n} z^{n}}e^{-\sum_{n=1}^{\infty}
{1\over n}b_{-n}w^{n}}\\  
&e&^{-\sum_{n=1}^{\infty}{q^{-2n}\over n(1+q^{2n})}b_{n}z^{-n}}  
e^{\sum_{n=1}^{\infty}{1\over  
n}b_{n}w^{-n}} \\ 
&\times& e^{-\al/2}w^{-\partial}(-q^3z)^{(\partial+i)/2}.  
\ea  

Let us now construct the VOS from the above VO. 
It is possible to  
consider the most general VOS, which is obtained from 
an arbitrary number of   
successive applications of the modes of $\phi^{i,-}(z)$ 
on the bosonic 
vacuum state. This, however,  is very complicated. 
We will instead
focus on the case where the VOS are obtained from the action 
of at most two $\phi^{i,-}(z)$ on the bosonic vacuum state, 
and where the 
dual VOS are obtained from the action of at most two 
$\phi^{i,+}(z)$ 
on the dual bosonic vacuum state. We believe that 
these cases embody 
many of the qualitative features of the general case. Before 
constructing the VOS, here are some operator product 
expansions that 
we shall be using:  
\bac  
\phi^{1-i,-}(z)\phi^{i,-}(w)&=&(-q^{3}z)^{1/2}{(q^{2}wz^{-1}; 
q^{4})_{\infty}\over  
(q^{4}wz^{-1};q^{4})_{\infty}}:\phi^{1-i,-}(z)\phi^{i,-}(w):,\\  
E^-(z)\phi^{i,-}(w)&=&{1\over z(1-q^{4}wz^{-1})}:
E^-(z)\phi^{i,-}(w):,\\  
\phi^{i,-}(z)E^-(w)&=&-{1\over zq^{3}
(1-q^{-2}wz^{-1})}:\phi^{i,-}(z) 
E^-(w):,\\  
E^-(z)E^-(w)&=&z^2(1-wz^{-1})(1-q^2wz^{-1}):E^-(z)E^-(w):, 
\label{prelim}  
\ea  
where in $(1-z)^{-1}$, $|z|<1$ is meant, and where 
\be
(x,q)_\infty=\prod_{n=0}^{\infty}(1-xq^n).
\ee  
Let us recall that the vacuum states $1$ and 
$e^{\al/2}$ are identified 
with the highest weight vectors of the modules 
$V(\Lambda_0)$ and $V(\Lambda_1)$ 
respectively. Consequently, only $\phi^{0,\pm}(z)$ act 
on $1$, and 
only $\phi^{1,\pm}(z)$ act on $e^{\al/2}$. 
Since $1$ is annihilated 
by $\partial$ and $b_n$, $n>0$, and since $e^{\al/2}
.1=e^{\al/2}$, we find that $\phi^{0,-}(z)$ acts on 
$1$ according to:  
\bac 
\phi^{0,-}(z).1&=& \\ 
&e&^{\sum_{n=1}^{\infty}{q^{4n}\over n(1+q^{2n})} 
b_{-n}z^{n}}e^{\al/2}
=\prod_{n=1}^{\infty}\sum_{m_n=0}^{\infty}{b_{-n}^{m_n}
q^{4nm_n}z^{nm_n}\over  
n^{m_n}(1+q^{2n})^{m_n}m_n!}e^{\al/2} \\
&=&\sum_{\la\in {\cal P}}
{b_{-\la}(q^4z)^{|\la|}\over {z}_{\la}(q)}e^{\al/2}
=\sum_{n\geq 0}Z_n(q^4z)^ne^{\al/2},  
\label{baba}
\ea  
The sum in (\ref{baba}) is over the set of all 
partitions ${\cal P}$, and to 
the partition $\la=(\la_1,\dots,\la_s)=  
(1^{m_1},2^{m_2},\dots,k^{m_k})$ we have associated the symbols  
\bac  
b_{-\la}&=&b_{-\la_1}b_{-\la_2}\dots b_{-\la_s},\\  
{z}_{\la}(q)&=&\prod_{i=1}^k i^{m_i}(1+q^{2i})^{m_i}m_i!,\\  
Z_n&=&\sum_{|\la|=n}{b_{-\la}\over {z}_{\la}(q)}.  
\label{Qp}\ea  
Here 
\be
|\lambda|=\sum_{i=1}^{s}\lambda_i=\sum_{i=1}^{k}{m_i}
\ee
is the weight of the partition $\lambda$.

Let $p_\lambda$ denote the power sum symmetric functions: 
$p_\lambda = p_{\lambda_1}p_{\lambda_2} \dots ,$ with $p_i = 
\sum_k x_k^i$.  Upon using the scalar product
\begin{equation}
< p_\lambda, p_\mu >_q 
= z_\lambda(q) \delta_{\lambda, \mu},
\end{equation}
the following identification can be made
\begin{equation}
b_{-\lambda} = p_\lambda
\end{equation}
As a result, the $Z_n$ become homogeneous symmetric functions 
and are recognized to be particular one-row Macdonald 
functions.  In general, the symmetric functions 
$P_\lambda(x;s,t)$ of Macdonald with two parameters 
$s$ and $t$, are orthogonal with respect to the 
scalar product
\begin{equation}
< p_\lambda, p_\mu >_{s,t}
= z_\lambda (s,t) \delta_{\lambda, \mu},
\end{equation}
where 
\begin{equation}
z_\lambda(s,t)
= \prod_{i=1}^k i^{m_i} \biggl( \frac{1-s^i}{1-t^i} 
\biggr)^{m_i} m_i!\, .
\end{equation}
The one-row functions $P_{(n)}(x;s,t)$ are proportional [4] to 
$\sum_{|\lambda| = n} z_\lambda(s,t)^{-1}p_\lambda(x)$.  
If we set $s = t^2 = q^4$, we see in fact that $< \ \ , 
\ \ >_{q^4,q^2} = < \ \ , \ \ >_q$ and that 
$z_\lambda(q^4,q^2) = z_\lambda (q)$.  
We thus observe that the $Z_n$ are one-row $q$-zonal   functions.

Let us now construct the one-row VOS. They are obtained 
from the modes of $\phi^{-,0}(z)$ as follows:  
\bac  
\phi^{0,-}_{-n}.1&=&{1\over 2\pi i}\oint dz z^{-n-1}
\phi^{0,-}(z).1\\
&=& {1\over 2\pi i}\oint dz z^{-n-1}\sum_{m\geq
0}Z_m(q^{4}z)^me^{\al/2}=  
q^{4n}Z_ne^{\al/2},  
\label{1VOS}
\ea  
with $n>0$. These states belong to $V(\Lambda_1)$. Next, we 
construct 
the most general two-row VOS from the modes of $\phi^{0,-}(z)$. 
This 
time, these  VOS lie in $V(\Lambda_0)$. We need  
\bac  
:\phi^{1,-}(z)\phi^{0,-}(w): .1&=& (-q^3z)^{1/2}   
e^{\sum_{n=1}^{\infty}{q^{4n}\over n(1+q^{2n})}b_{-n}z^{n}}    
e^{\sum_{n=1}^{\infty}{q^{4n}\over n(1+q^{2n})}b_{-n}w^{n}}
e^\al\\  
&=&  
(-q^3z)^{1/2}\sum_{m_1,m_2\geq
0}q^{4(m_1+m_2)}Z_{m_1}Z_{m_2}z^{m_1}  
w^{m_2}e^\al.  
\ea  
Using this relation and (\ref{prelim}), we find
\bac
\lefteqn{ \phi^{1,-}_{-r}\phi^{0,-}_{-s} .1 ={1\over (2\pi i)^2}
\oint   
dz dw z^{-r-1} w^{-s-1}\phi^{1,-}(z)\phi^{0,-}(w).1}\\  
&&\qquad ={1\over (2\pi i)^2}\oint   
dz dw z^{-r-1} w^{-s-1}(-q^{3}z)^{1/2}{(q^{2}wz^{-1};
q^{4})_{\infty}\over  
(q^{4}wz^{-1};q^{4})_{\infty}}:\phi^{1,-}(z)\phi^{0,-}(w):.1\\  
&&\qquad=-{q^3\over (2\pi i)^2}\oint   
dz dw {(q^{2}wz^{-1};q^{4})_{\infty}\over  
(q^{4}wz^{-1};q^{4})_{\infty}}  
\sum_{m_1,m_2\geq 0}q^{4(m_1+m_2)}Z_{m_1}Z_{m_2}z^{m_1-r}  
w^{m_2-s-1}e^\al\\  
&&\qquad=-{q^3\over (2\pi i)^2}\oint   
dz  dw\sum_{m_1,m_2,n\geq 0}q^{4(m_1+m_2)}C_{n}Z_{m_1}Z_{m_2}
z^{m_1-r-n}  
w^{m_2-s-1+n}e^\al\\  
&&\qquad=-{q^3\over 2\pi i}\oint   
dz \sum_{m_1,n\geq 0}q^{4(m_1+s-n)}C_{n}Z_{m_1}Z_{s-n}
z^{m_1-r-n} e^\al\\  
&&\qquad=-q^3 \sum_{n\geq 0}q^{4(r+s-1)}C_{n}Z_{n+r-1}Z_{s-n}  
e^\al\\  
&&\qquad =-q^{4(r+s)-1} \sum_{n\geq 0}C_{n}(R_{12})^n
Z_{r-1}Z_{s}  
e^\al=-q^{4(r+s)-1} {(q^2R_{12};q^4)_{\infty}\over   
(q^4R_{12};q^4)_{\infty}}Z_{r-1}Z_{s}e^\al,  
\label{ff} 
\ea
where $R_{12}$ is the raising operator acting on 
$Z_{r}Z_{s}$ according to
\be  
R_{12}(Z_{r}Z_s)=Z_{r+1}Z_{s-1},  
\label{ro} 
\ee  
and the coefficients $C_n$, $n\geq 0$ are such that  
\be  
\sum_{n\geq 0}C_n x^n= {(q^2x;q^4)_{\infty}\over   
(q^4x;q^4)_{\infty}}=e^{-\sum_{k\geq 1}{q^{2k}x^k\over 
k(1+q^{2k})}}=  
\sum_{\la\in {\cal P}}{(-1)^{\ell(\la)}(q^{2}x)^{|\la|} \over 
{z}_{\la}(q)},  
\ee  
that is,
\be  
C_n=\sum_{|\la|=n}{(-1)^{\ell(\la)}q^{2|\la|}\over 
{z}_{\la}(q)}.  
\label{CC}\ee  
Let us now consider the construction of a single-row dual 
VOS from the  
action of the conjugate VO $\phi^{1,+}(z^{*-1})$ on the dual 
bosonic 
vacuum 1. For this purpose let us remark that, unlike the 
situation in the 
classical case (i.e., $q=1$), the dual VOS cannot simply be 
obtained  
from the VOS through the usual conjugation $b_n^*=b_{-n}$ 
because
this leads to a dual Fock space that cannot be identified 
with any 
of the dual modules of $V(\Lambda_i)$. For this reason, 
we construct 
the dual VOS by operating with $\phi^{1,+}(z^{*-1})$. 
The action of the 
modes $b_n$, $\partial$, and $e^\al$ on the dual bosonic 
vacuum are 
given by   
\bac  
1. b_{-n}&=&1.\partial =0,\quad n>0,\\  
1. e^{-\al}&=&e^{-\al},  
\ea  
and $1. b_n$, $n>0$ are non-zero states. Due to the integral 
form of 
this VO, the following new technical problem arises: the 
contour in 
the definition of $\phi^{1,+}(z^{*-1})$ (\ref{vo}) winds 
around the 
pole $w=0$ and it is not easy to compute the residue at this 
pole in 
a closed form. Therefore, we make a change of variable to 
remove this 
pole when $\phi^{1,+}(z^{*-1})$ acts on the dual bosonic vacuum
1.   
The appropriate change of variable is simply $w=\xi^{-1}$. 
Let us 
also set $z^{*-1}=\eta$, then we  obtain  
\bac  
1.:\phi^{1,-}(\eta)E^-(\xi^{-1}):&=&  
1. e^{-\sum_{n=1}^{\infty}{q^{-2n}\over 
n(1+q^{2n})}b_{n}\eta^{-n}}  
e^{\sum_{n=1}^{\infty}{1\over  
n}b_{n}\xi^{n}}e^{-\al/2}
\xi^{\partial}(-q^3\eta)^{(\partial+1)/2}\\  
&=&  
e^{-\al/2}e^{-\sum_{n=1}^{\infty}{q^{-2n}
\over n(1+q^{2n})}b_{n}\eta^{-n}}  
e^{\sum_{n=1}^{\infty}{1\over  
n}b_{n}\xi^{n}}\xi(-q^3\eta)\\  
1.:\phi^{1,-}(\eta)E^-(\eta q^2):&=&-qe^{-\al/2}  
e^{\sum_{n=1}^{\infty}{1\over n(1+q^{2n})}b_{n}\eta^{-n}}.  
\ea  
From these relations we arrive at 
\bac  
1. \phi^{1,+}(\eta)&=&  
-{\q\over 2\pi i}\oint_{|\eta^{-1}q^{-2}|<|\xi|<|\eta^{-1} 
q^{-4}|} d\xi {  
1. :\phi^{1,-}(\eta)  
E^-(\xi^{-1}):\over   
\xi\eta^2q^{6}(\xi-\eta^{-1}q^{-2})(\xi-\eta^{-1} q^{-4})},\\  
&=&-q^{-1}1.:\phi^{1,-}(\eta)E^-(\eta q^2):=e^{-\al/2}  
e^{\sum_{n=1}^{\infty}{1\over n(1+q^{2n})}b_{n}\eta^{-n}}\\  
&=&e^{-\al/2}\sum_{\la\in {\cal P}}{b_{\la}\over {z}_{\la}(q)}
\eta^{-|\la|}=  
e^{-\al/2}\sum_{n\geq 0}Z^*_n\eta^{-n},  
\label{1DVOS}\ea  
where  
\bac  
b_{\la}&=&b_{\la_1}b_{\la_2}\dots,\\  
Z^*_n&=&\sum_{|\la|=n}{b_{\la}\over {z}_{\la}(q)}.  
\ea  
Keeping with (2.15), we make the identification  
\bac 
b_n&=&n(1+q^{2n}){\partial\over \partial p_n}\\ 
&=&D(p_n),\quad n>0. 
\ea 
Here the first equality is consistent with the Heisenberg 
algebra  
(\ref{HA}) and the second one defines the adjoint operator  
$D$ \cite{JaYu93} with respect to the scalar product (2.14) 
on the space of  symmetric functions.  We thus have  
\be 
Z^*_\lambda=D(Z_\lambda). 
\ee 
From this, we obtain the following scalar product  
\be  
<Z_n,Z_m>_q=\delta_{n,m}\sum_{|\la|=n}{1\over {z}_{\la}(q)}.  
\label{sca}\ee  
The one-row dual VOS are hence found to be  
\be  
1. \phi^{1,+}_n={1\over 2\pi i}\oint d\eta \eta^{n-1}1.
\phi^{1,+}(\eta)=   
{e^{-\al/2}\over 2\pi i}\oint d\eta \eta^{n-1}\sum_{m\geq 0}
Z^*_m\eta^{-m}  
=e^{-\al/2}Z^*_n.  
\label{re}
\ee  
As a check, we now show that the scalar products of 1-row VOS 
are consistent with the q-KZ equation \cite{FR}. 
Indeed, from (\ref{1VOS}), 
(\ref{sca}) and 
(\ref{re}) we get the following matrix element:
\be
1. \phi_m^{1,+}\phi_{-n}^{0,-} .1 =q^{4n}<1,Z^*_mZ_n>_{q}=
\delta_{n,m}q^{4n}\sum_{|\la|=n}{1\over {z}_{\la}(q)},
\ee
which in turn leads to 
\bac
1.\phi^{1,+}(z)\phi^{0,-}(w) .1 &=&
\sum_{n,m\in  \Z}z^{-m}w^{n}
1.\phi_m^{1,+}\phi_{-n}^{0,-} .1\\
& =&
\sum_{n\in  \Z}\sum_{|\la|=n}
z^{-n}w^{n}
{q^{4n}\over {z}_{\la}(q)}
=
{(q^6wz^{-1};q^4)_{\infty}\over   
(q^4wz^{-1};q^4)_{\infty}},
\ea
but this is precisely the matrix element of the VO obtained 
through 
solving the q-KZ in \cite{DFJMN}.
 
Let us now construct the most general two-row dual VOS from the  
successive actions of two modes, one from $\phi^{1,+}(z^{*-1})$  
and the other from  $\phi^{0,+}(w^{*-1})$, on the dual bosonic 
vacuum. Let $\eta=z^{*-1}$ and $\theta=w^{*-1}$, we find that
\be 
1. \phi^{1,+}(\eta)\phi^{0,+}(\theta)=q^{-4} 
{(q^2\theta\eta^{-1};q^4)_{\infty}\over 
(q^4\theta\eta^{-1};q^4)_{\infty}} 
e^{-\alpha}\sum_{n,m\geq 0}Z^*_nZ^*_m \eta^{-n}\theta^{-m-1}. 
\ee 
Thus 
\bac  
1. \phi^{1,+}_{r}\phi^{0,+}_{s}&=& 
{1\over (2\pi i)^2}\oint   
d\eta d\theta \eta^{r-1} \theta^{s-1}1.\phi^{1,+}(\eta)
\phi^{0,+}(\theta)\\  
&=&{q^{-4}e^{-\alpha}\over (2\pi i)^2}\oint   
d\eta d\theta \eta^{r-1} \theta^{s-1}{(q^{2}\theta  
\eta^{-1};q^{4})_{\infty}\over  
(q^{4}\theta \eta^{-1};q^{4})_{\infty}} 
\sum_{n,m\geq 0}Z^*_nZ^*_m \eta^{-n}\theta^{-m-1}\\  
&=&{q^{-4}e^{-\alpha}\over (2\pi i)^2}\oint   
d\eta d\theta  
\sum_{n,m,k\geq 0}C_kZ^*_nZ^*_m \eta^{r-n-k-1}
\theta^{s-m+k-2}\\  
&=&{q^{-4}e^{-\alpha}\over 2\pi i}\oint   
d\eta   
\sum_{n,k\geq 0}C_kZ^*_nZ^*_{s+k-1} \eta^{r-n-k-1}\\  
&=&q^{-4}e^{-\alpha}   
\sum_{k\geq 0}C_kZ^*_{r-k}Z^*_{s+k-1}= 
q^{-4}e^{-\alpha}   
\sum_{k\geq 0}C_k(\~R_{12})^kZ^*_{r}Z^*_{s-1} \\  
&=&q^{-4}e^{-\alpha}   
{(q^2\~R_{12};q^4)_{\infty}\over   
(q^4\~R_{12};q^4)_{\infty}}Z^*_{r}Z^*_{s-1}, 
\label{2DVOS}\ea 
where the coefficients $C_n$ are the same as those given in 
(\ref{CC}) and $\~R_{12}$ is the lowering operator defined by
\be
\~R_{12}Z_{r}Z_{s}=Z_{r-1}Z_{s+1}.
\ee
 
We know through the relations (\ref{1VOS}) and (\ref{1DVOS}) 
that the VOS obtained from the action of a single mode of 
the vertex operators on the vacuum are proportional to the 
one-row $q$-zonal  functions. This is no longer true in the case 
of VOS obtained from the action of two modes, that is, they 
are not proportional to two-row $q$-zonal  functions. 
The two are however   related and we now display this 
relation. For this 
purpose, we use the following relation obtained in
 \cite{JJ},  which expresses the two-row $q$-zonal    
functions $Z_{n-1,m}$ as the product of  two 
one-row $q$-zonal  functions: 
\be 
Z_{n-1,m}=\>  _{2}\phi_1(q^2,q^{4(n-m)};q^{2+4(n-m)};q^4,q^2 R) 
(1-R)Z_{n-1}Z_m. 
\label{JJ}
\ee 
We have set  $R=R_{12}$ in (2.37) and as usual 
$_{2}\phi_1(a,b;c;q,z)$ denotes  the $q$-hypergeometric 
function  
\be 
_{2}\phi_1(a,b;c;q,z)=\sum_{n\geq 0}{(a;q)_n(b;q)_n\over  
(c;q)_n}{z^n\over (q;q)_n}, 
\ee 
with 
\be 
(a,q)_n=(1-a)\dots (1-aq^{n-1}). 
\ee  
Combining the relations (\ref{JJ}) and (\ref{ff}), we get  
\be 
Z_{n-1,m}=-q^{-4(n+m)+1} \> _{2}
\phi_1(q^2,q^{4(n-m)};q^{2+4(n-m)};q^4,q^2 R) 
{(R;q^4)_{\infty}\over (q^2R;q^4)_{\infty}} 
e^{-\alpha}\phi^{1,-}_{-n}\phi^{0,-}_{-m}.1, 
\ee 
which expresses the two-row $q$-zonal  
functions in terms of the VOS. 
To also relate the dual VOS (\ref{2DVOS})   
to the two-row $q$-zonal  functions, 
we use the relation 
\be 
(Z_{n-1,m})^*=\> _{2}\phi_1(q^2,
q^{4(n-m)};q^{2+4(n-m)};q^4,q^2 \~R) 
(1-\~R)(Z_{m})^*(Z_{n-1})^*, 
\ee 
which is deduced from (\ref{JJ}). Then, (\ref{2DVOS}) leads to 
\be 
(Z_{n-1,m})^*=-q^{4}\> _{2}\phi_1(q^2,q^{4(n-m)};q^{2+4(n-m)};
q^4,q^2 \~R) 
{(\~R;q^4)_{\infty}\over (q^2\~R;q^4)_{\infty}} 
1.\phi^{1,+}_{m}\phi^{0,+}_{n}e^{\alpha}, 
\ee 
where $\~R=\~R_{12}$.
 
\section{Discussion and conclusions} 
 
Let us compare our results with those of Jing  who 
considers in \cite{Jing1} the following vertex operators: 
\be  
X(z)=e^{\sum_{n=1}^{\infty}{1-t^{n}\over n(1-s^n)}c_{-n}z^{n}}  
e^{-\sum_{n=1}^{\infty}{1-t^{n}\over  
n(1-s^n)}c_{n}z^{-n}}=\sum_{n\in\Z}X_{n}  
z^{-n},  
\label{Ji}  
\ee  
where $s$ and $t$ are deformation parameters and where 
the bosonic modes 
$c_n$ satisfy the Heisenberg algebra defining relation: 
\be  
{[c_n, c_m]}=n{1-s^{|n|}\over 1-t^{|n|}}\delta_{n+m,0}.  
\ee  
The modes $X_n$ are expressed in terms  of $X(z)$ as  follows 
\be  
X_n={1\over 2\pi i}\oint dz z^{n-1}X(z).  
\ee  
Jing   defines the conjugate vertex operator $X^*(z)$ by  
\be  
X^*(z)=e^{-\sum_{n=1}^{\infty}{1-t^{n}\over n(1-s^n)}
c_{-n}z^{n}}  
e^{\sum_{n=1}^{\infty}{1-t^{n}\over n(1-s^n)}c_{n}z^{-n}}= 
\sum_{n\in\Z}X^*_{n}  
z^{n},  
\ee  
and he also defines the adjoint of $c_n$ in the Fock space by  
$c_n^*=c_{-n}$. It can thus easily be checked that   
\be  
(X(z))^*=X^*(z^{*-1}),  
\ee  
from where it follows that  
\be  
(X_n)^*=X^*_n.  
\ee  
Thus, $X^*_n$ is the conjugate of $X_n$ and for this reason we  
call $X^*(z)$  the conjugate VO of $X(z)$.   
 
Let us note  that the Heisenberg algebra elements $c_n$  
reduce to the elements $b_n$ when $s=t^2=q^4$. However, there 
are three  differences between the intertwining VO and those 
considered by Jing. Firstly, the conjugate VO $X^*(z)$ has the 
same simple form as the direct VO $X(z)$, whereas 
$\phi^{1,+}(z)$ which is the conjugate VO of $\phi^{0,-}(z)$ 
does not have the same simple exponential form as that of 
$\phi^{0,-}(z)$. Secondly, Jing's deformation is symmetric in 
the following sense: if   
\be  
{[c_n,c_m]}=nf(s^{|n|}, t^{|n|})\delta_{n+m,0},  
\ee  
with 
\be
f(s^{|n|}, t^{|n|})={1-s^{|n|}\over 1-t^{|n|}},
\ee
then the deformation of the VO takes the form  
\be  
X(z)=  
e^{\sum_{n=1}^{\infty}{1\over nf(s^{n},t^{n})}c_{-n}z^{n}}  
e^{-\sum_{n=1}^{\infty}{1\over nf(s^{n}, t^{n})}c_{n}z^{-n}}.  
\ee   
This is of course not the case for the intertwining VO. Finally, 
$X(z)$ and $X^*(z)$ do not depend on the zero mode $\partial$ 
and its conjugate $\al$ whereas the intertwining VO's do. For 
these reasons Jing's VOS and our VOS are related differently   
to two-row $q$-zonal  functions \cite{Jing1}.

One important open question is to understand the relation 
between 
the matrix elements of the VOS (rather than the VOS themselves) 
and   symmetric functions. It is known in fact that 
these matrix elements satisfy the q-KZ equation \cite{FR} and 
  it would be quite useful to find a  relation between this 
equation 
and the symmetric functions. Another interesting question is to
understand the relation between the intertwining vertex 
operators and the creation operators recently 
introduced in [10] which allow to construct from the 
ground state of the Calogero-Sutherland model the 
Jack functions associated to any partition.

A straightforward extension of this work is to consider the 
other 
type of VO $\psi^{i,\pm}(z)$, referred to as type-II in 
\cite{DFJMN}  that intertwine the $\uq$ modules as follows:

\be  
v_{+}\otimes \psi^{i,+}(z)+v_{-}\otimes \psi^{i,-}(z):\quad  
V(\Lambda_i)\ra V(z)\otimes V(\Lambda_{i-1}).  
\ee  

Furthermore, it would be interesting to carry this work over 
to the VO which intertwine the higher level modules of $\uq$ 
\cite{BW,B}. Finally, and most importantly, extending this 
work beyond 
the two-row $q$-zonal  functions, is certainly 
an interesting open problem which is also pending in the 
framework of
 Jing. 
 
\section{\bf Acknowlegement} 

A.H.B. is supported by an NSERC postdoctoral fellowship. 
L.V. is supported by NSERC (Canada) and FCAR (Qu\'ebec). 
We are very grateful to O. Foda for very stimulating  
discussions and
for reading this paper. 
We wish to thank T. Baker, P. Jarvis and N. Jing for drawing 
our attention to many of their interesting papers. 
   
\end{document}